\documentclass[12pt]{amsart}
\usepackage{amsmath,amsfonts,amssymb,graphicx,hyperref, enumitem,amsthm,mathabx,amsrefs}
\hypersetup{colorlinks=true,citecolor=blue,linkcolor=blue,urlcolor=blue,pdfstartview=FitH}

\setlist[enumerate,1]{label={\arabic*.}}

\newcommand{\Z}{\mathbb{Z}}
\newcommand{\R}{\mathbb{R}}
\newcommand{\N}{\mathbb{N}}
\newcommand{\C}{\mathbb{C}}

\newcommand{\paren}[1]{\ensuremath{\left( #1 \right)}}
\newcommand{\set}[1]{\ensuremath{\left\{ #1 \right\}}}
\newcommand{\abs}[1]{\ensuremath{\left| #1 \right|}}
\newcommand{\norm}[1]{\ensuremath{\left\| #1 \right\|}}
\newcommand{\setdiv}{\,\middle|\,}
\newcommand{\summod}[1]{\ensuremath{\,(\mathrm{mod}\,#1)}}
\newcommand{\e}[1]{e\paren{#1}}
\DeclareMathOperator*{\res}{res}
\DeclareMathOperator{\diag}{diag}
\renewcommand{\Re}{{\mathop{\mathgroup\symoperators Re}}}
\renewcommand{\Im}{{\mathop{\mathgroup\symoperators Im}}}
\newcommand{\sgn}{{\mathop{\mathgroup\symoperators \,sgn}}}
\newcommand{\wbar}[1]{\overline{#1}}
\newcommand{\wtilde}[1]{\widetilde{#1}}

\newcommand{\Ad}{\operatorname{Ad}}

\newcommand{\specmu}{\mathbf{spec}}
\newcommand{\cosmu}{\mathbf{cos}}
\newcommand{\sinmu}{\mathbf{sin}}

\theoremstyle{plain} 
\newtheorem{thm}{Theorem}

\newtheorem{prop}[thm]{Proposition}

\title{The arithmetic Kuznetsov formula on $GL(3)$, I:\\ The Whittaker case.}
\author{Jack Buttcane}
\date{29 May 2018}
\address{Mathematics Department, 244 Mathematics Building, Buffalo, NY 14260, USA}
\email{buttcane@buffalo.edu}
\thanks{During the time of this research, the author was supported by NSF grant DMS-1601919.}

\begin{document}

\begin{abstract}
The original formulae of Kuznetsov for $SL(2,\Z)$ allowed one to study either a spectral average via Kloosterman sums or to study an average of Kloosterman sums via a spectral interpretation.
In previous papers, we have developed the spectral Kuznetsov formulae at the minimal weights for $SL(3,\Z)$, and in these formulae, the big-cell Kloosterman sums occur with weight functions attached to four different integral kernels, according to the choice of signs of the indices.
These correspond to the $J$- and $K$-Bessel functions in the case of $GL(2)$.
In this paper, we demonstrate a linear combination of the spherical and weight-one $SL(3,\Z)$ Kuznetsov formulae that isolates one particular integral kernel, which is the spherical $GL(3)$ Whittaker function.
Using the known inversion formula of Wallach, we give the first arithmetic Kuznetsov formula for $SL(3,\Z)$ and use it to study smooth averages and the Kloosterman zeta function attached to this particular choice of signs.
\end{abstract}


\maketitle

\section{Introduction}

The meromorphic continuation of the Poincar\'e series for $SL(2,\Z)$ was first conducted by Selberg \cite{Sel04} for the purposes of bounding Fourier coefficients of automorphic forms.
The Fourier coefficients of the $SL(2,\Z)$ Poincar\'e series can be expressed in terms of the Kloosterman zeta function
\[ \widetilde{Z}_2(s) := \sum_{c = 1}^\infty \frac{S(m,n;c)}{c^{1+s}}, \]
where
\[ S(m,n;c) = \sum_{\substack{x\summod{c}\\x\bar{x}\equiv 1 \summod{c}}} \e{\frac{mx+n\bar{x}}{c}}, \qquad \e{t} = e^{2\pi i t} \]
is the classical Kloosterman sum.
Weil's bound $S(m,n;c)\ll_{m,n} c^{\frac{1}{2}+\epsilon}$ implies the Kloosterman zeta function converges absolutely on $\Re(s) > \frac{1}{2}$, and Selberg's study of Poincar\'e series implies it has meromorphic continuation to all $s$ with poles at integral shifts of the spectral parameters of the $SL(2,\Z)$ Maass cusp forms, as well as certain poles arising from the continuous spectrum.

Somewhat later, Kuznetsov \cite{Kuz01} used their spectral expansion to resolve a conjecture of Linnik \cite{L01}, itself motivated by problems in additive number theory:
\begin{thm}[Kuznetsov]
	For $X > 1$ a large parameter and $m,n\in\Z$, $mn > 0$,
\begin{align*}
	\sum_{c \le X} \frac{S(m,n;c)}{c} \ll_{m,n} X^{\frac{1}{6}} \paren{\log X}^{\frac{1}{3}},
\end{align*}
\end{thm}
For contrast, Weil's square-root cancellation bound would give $X^{\frac{1}{2}+\epsilon}$, so we are seeing cancellation between terms in the sum over the moduli $c$.
Away from $SL(2,\Z)$ where the Selberg eigenvalue conjecture is known, generalizations of Kuznetsov's result to congruence subgroups typically involve an additional term $X^\theta$ where $0 < \theta < \frac{1}{2}$ comes from bounds on the spectral parameters of the relevant Maass cusp forms.

Kuznetsov's proof came from the construction of two formulae that strongly tie moduli sums of Kloosterman sums to the study of automorphic forms.
The first, which we call the spectral Kuznetsov formula, expresses certain averages of a test function over $SL(2,\Z)$ automorphic forms as sums of Kloosterman sums, weighted by integral transforms of the test function.
The second, which we call the arithmetic Kuznetsov formula, operates in reverse, expressing averages of Kloosterman sums weighted by a test function as sums over automorphic forms weighted by integral transforms of the test function.
The first direction may be thought of as a generalization of Petersson's trace formula for modular forms to the $SL(2,\Z)$ Maass forms.

The first breakthrough in the study of Poincar\'e series on higher-rank groups came in the paper \cite{BFG}.
Among other things, that paper proves in detail the meromorphic continuation of a spherical $SL(3,\Z)$ Poincar\'e series.
The meromorphic continuation of the corresponding Kloosterman zeta function is also claimed \cite{BFG}*{eq. (1.20)}, but the proof was never provided.

In the papers \cite{SpectralKuz, WeylI, WeylII}, the author completed the generalization of the spectral Kuznetsov formula to $SL(3,\Z)$ and gave the non-spherical analogs, i.e. Kuznetsov-type trace formulae for automorphic forms which are non-trivial on $SO(3,\R)$.
In this paper, we turn our attention to building the arithmetic Kuznetsov formulae on $SL(3,\Z)$.
These will consider sums of the long-element $SL(3, \Z)$ Kloosterman sum.
We denote that sum by $S_{w_l}(\psi_m,\psi_n,c)$, where $m,n\in\Z^2$ are the indices and $c\in\N^2$ are the moduli.
The definition of this exponential sum is somewhat unpleasant and not entirely relevant to the discussion, so we postpone that until Section \ref{sect:Results}, below.

We have square-root cancellation bounds originally due to Stevens \cite{Stevens} of the form
\begin{align}
\label{eq:StevensBd}
	S_{w_l}(\psi_m,\psi_n,c) \ll_{m,n} \sqrt{(c_1,c_2)} (c_1 c_2)^{\frac{1}{2}+\epsilon},
\end{align}
though \cite{DabFish} has improved the $\sqrt{(c_1,c_2)}$ term (and much more explicit bounds are given in \cite{Me01}*{Theorem 2} and \cite{LargeSieve}*{eq. (2.10)}).

The author's thesis \cite{Me01} provides a first attempt at an arithmetic Kuznetsov formula using a first-term inversion formula on the spherical Kuznetsov formula.
\begin{thm}
\label{thm:Me01}
Let $X_1, X_2 > 1$ be large parameters.
Suppose $m,n\in \Z^2$ such that $m_1 n_2 m_2 n_1 \ne 0$, and $f$ smooth and compactly supported on $(\R^+)^2$, then
\begin{align*}
	\sum_{\varepsilon \in \set{\pm1}^2} \sum_{c_1,c_2\in\N} \frac{S_{w_l}(\psi_m,\psi_{\varepsilon n},c)}{c_1 c_2} f\paren{\tfrac{X_1 c_2}{c_1^2}, \tfrac{X_2 c_1}{c_2^2}} \ll_{m,n,f,\epsilon}& (X_1 X_2)^{\theta+\epsilon} + X_1^{\frac{1}{2}+\epsilon} + X_2^{\frac{1}{2}+\epsilon},
\end{align*}
where $\theta$ is any bound toward the Ramanujan-Selberg conjecture for spherical $SL(3,\Z)$ Maass forms.
\end{thm}
The Kim-Sarnak bound \cite{KS01} shows $\theta=\frac{5}{14}$ is acceptable.
The square-root cancellation bound \eqref{eq:StevensBd} implies the sum is at most $(X_1 X_2)^{\frac{1}{2}+\epsilon}$, so we are seeing cancellation between terms of the $c$ sum, provided the moduli are not too far apart.
That is, the terms $X_1^{\frac{1}{2}+\epsilon} + X_2^{\frac{1}{2}+\epsilon}$ become comparable to the square-root cancellation bound when, say, $c_1 < c_2^{1/2+\epsilon}$.

The method of \cite{Me01} can be regarded, in a much simplified form, as writing a Fourier coefficient of a spherical $SL(3,\Z)$ Poincar\'e series in the form
\[ \widetilde{Z}(s) + F(s_1+\tfrac{1}{3},s_2-\tfrac{1}{6}) + G(s_1-\tfrac{1}{6},s_2+\tfrac{1}{3}), \]
where $\widetilde{Z}(s)$, $F(s)$ and $G(s)$ are holomorphic on $\Re(s_1),\Re(s_2) > \frac{1}{6}$, and
\[ \widetilde{Z}(s) = \sum_{w \in W} C_w(m,n,s) \widetilde{Z}_w(s), \]
is the full Kloosterman zeta function with $W$ the Weyl group, $C_w(m,n,s)$ a quotient of gamma functions (and powers of $m$ and $n$) which varies with $w\in W$, and $\widetilde{Z}_w(s)$ the sign-independent Kloosterman zeta function at each Weyl element; in particular,
\begin{align}
\label{eq:SignIdepKZ}
	\widetilde{Z}_{w_l}(s) = \sum_{\varepsilon \in \set{\pm1}^2} \sum_{c_1,c_2\in\N} \frac{S_{w_l}(\psi_m,\psi_{\varepsilon n},c)}{c_1^{1+3s_1} c_2^{1+3s_2}}.
\end{align}
The extra terms $F$ and $G$ here are the direct cause of the terms $X_1^{\frac{1}{2}+\epsilon} + X_2^{\frac{1}{2}+\epsilon}$ in Theorem \ref{thm:Me01}; the strange combination of shifts is because we save factors of $\frac{\sqrt{c_2}}{c_1}$ and $\frac{\sqrt{c_1}}{c_2}$ instead of $\frac{1}{c_1}$ and $\frac{1}{c_2}$ directly.
The term $F(s_1+\frac{1}{3},s_2-\frac{1}{6})$ gives some analytic continuation in the variable $\tilde{s}_1 := 2s_1+s_2$ and $G(s_1-\frac{1}{6},s_2+\frac{1}{3})$ gives some analytic continuation in $\tilde{s}_2 := 2s_2+s_1$, but the sum of the two terms gives no analytic continuation in either variable.
It is important to note that the functions $F$ and $G$ are not sums of shifts of $\widetilde{Z}$; this is because the kernel functions vary with the signs $\varepsilon$, even though their first-term asymptotics do not.
In this way, it is not possible to use the methods of \cite{Me01} to obtain the meromorphic continuation of the long-element Kloosterman zeta functions.
We expect that the methods of \cite{BFG} would encounter similar difficulties.

Here we will separate out a particular choice of signs and apply a true inversion formula.
As we will see, it is necessary to have both the spherical and weight-one Kuznetsov formulae in order to accomplish this.

We treat the long-element $SL(3,\Z)$ Kloosterman sum $S_{w_l}(\psi_m,\psi_n,c)$ in the case \\ $m_1 n_2, m_2 n_1 > 0$.
When the moduli are coprime, the long-element Kloosterman sum factors into classical Kloosterman sums as
\[ S_{w_l}(\psi_m,\psi_n,c) = S(m_1,-n_2 c_2; c_1) S(m_2 c_1, -n_1; c_2) \]
so one might (counter-intuitively) term this as the $-,-$ case for the Kloosterman sums.

The arithmetic Kuznetsov formula attaches a spectral interpretation to smooth averages of Kloosterman sums.
We state the full formula in Theorem \ref{thm:ArithKuz}, below.
A simple application gives the following bound:
\begin{thm}
\label{thm:SmoothSums}
Let $X_1, X_2 > 0$ with $X_1 X_2 > 1$.
Suppose $m,n\in \Z^2$ such that $m_1 n_2, m_2 n_1 > 0$, and $f$ smooth and compactly supported on $(\R^+)^2$, then
\begin{align*}
	\sum_{c_1,c_2\in\N} \frac{S_{w_l}(\psi_m,\psi_n,c)}{c_1 c_2} f\paren{\tfrac{X_1 c_2}{c_1^2}, \tfrac{X_2 c_1}{c_2^2}} \ll_{m,n,f,\epsilon}& (X_1 X_2)^{\theta+\epsilon}+X_1^{-1-\epsilon}+X_2^{-1-\epsilon},
\end{align*}
where $\theta$ is any bound toward Ramanujan-Selberg for the spherical and weight-one $SL(3,\Z)$ Maass forms.
\end{thm}
We note below that the Kim-Sarnak result $\theta=\frac{5}{14}$ extends to these forms, as well.
Again, the square-root cancellation bound is $(X_1 X_2)^{\frac{1}{2}+\epsilon}$, so we are seeing cancellation between terms of the $c$ sum, even when the moduli are very far apart.
In fact, the terms $X_1^{-1-\epsilon}+X_2^{-1-\epsilon}$ only become comparable to the square-root cancellation bound when, say $c_1 < c_2^{1/5+\epsilon}$.
As is typical for smooth averages, the Ramanujan-Selberg conjecture $\theta=0$ implies the bound $(X_1 X_2)^\epsilon+X_1^{-1-\epsilon}+X_2^{-1-\epsilon}$.
Aside from the separation of signs, the bound in this theorem is stronger than that of Theorem \ref{thm:Me01} due to the use of the full inversion formula, Theorem \ref{thm:KontLeb}, which effectively removes the $X_1^{\frac{1}{2}+\epsilon} + X_2^{\frac{1}{2}+\epsilon}$ terms.
The terms $X_1^{-1-\epsilon}+X_2^{-1-\epsilon}$ derive from a particularly vexing group of Maass cusp forms in the neighborhood of the self-dual forms, which will haunt us throughout the paper.

We also obtain the meromorphic continuation of the simplest weighted Kloosterman zeta function for this choice of signs:
\begin{thm}
\label{thm:KloosZeta}
	Let $m,n\in \Z^2$ with $m_1 n_2, m_2 n_1 > 0$, and for $s\in\C^2$, define
	\begin{align}
	\label{eq:tildesdef}
		\tilde{s}=&(2s_1+s_2,s_1+2s_2),
	\end{align}
	then the Kloosterman zeta function
	\begin{align*}
		Z_{m,n}^*(s) :=& \sum_{c_1,c_2\in\N} \frac{S_{w_l}(\psi_m,\psi_n,c)}{c_1 c_2} f_s\paren{\tfrac{2\sqrt{m_1 n_2 c_2}}{c_1}, \tfrac{2\sqrt{m_2 n_1 c_1}}{c_2}}, \\
		f_s(y) =& y_1^{2\tilde{s}_1} y_2^{2\tilde{s}_1} \exp\paren{-\pi^2(y_1^2+y_2^2)},
	\end{align*}
	initially convergent on $\Re(s_1), \Re(s_2) > \frac{1}{2}$, has meromorphic continuation to all of $s \in \C^2$ with poles whenever
	\begin{align}
	\label{eq:KloosZetaPoles}
	\begin{aligned}
		\tilde{s}_1=&-\mu_i-\ell, &\text{or}&&\tilde{s}_2 =& \mu_i-\ell,
	\end{aligned}
	\end{align}
	 with $\mu=\mu_\varphi$ for some cusp form $\varphi$ of weight at most one, $i=1,2,3$ and $\ell \in \N_0=\N\cup\set{0}$.
\end{thm}
The above theorem does not list the poles coming from the Eisenstein series, of which there are many, as these are somewhat more complicated.
In Section \ref{sect:KloosZeta} we discuss the orders and residues of the cuspidal poles and briefly indicate how one may determine the poles of the Eisenstein series terms.

The weighted Kloosterman zeta function used in the theorem does not lend itself to proving the meromorphic continuation of the unweighted Kloosterman zeta function, due to the same complications described above.
It seems that even having the arithmetic Kuznetsov formula is not sufficient to obtain said continuation, and we discuss some reasons why this fails in Section \ref{sect:KloosZeta}.

The reason why our arithmetic Kuznetsov formula in the $m_1 n_2, m_2 n_1 > 0$ case involves only the spherical and weight-one forms is because that term is zero in the remaining $d \ge 2$ spectral Kuznetsov formulae \cite{WeylII}, for algebraic reasons.
This behavior is very similar to the case of $SL(2,\Z)$, where the Whittaker function is the $K$-Bessel function.
It is not immediately clear how to isolate the remaining signs, but any arithmetic Kuznetsov formulae in those cases will certainly involve all of the spectral Kuznetsov formulae and will require constructing inversion formulae with discrete spectra.
The meromorphic continuation of the sign-independent Kloosterman zeta function \eqref{eq:SignIdepKZ}, if it exists, should then also involve the Maass cusp forms of every weight.

An alternative to the current approach to studying sums of Kloosterman sums is to define the Poincar\'e series directly in terms of the Bruhat decomposition (see \cite{CPS}); this bypasses the need for the inversion of integral transforms.
The construction of such Poincar\'e series on $SL(n,\Z)$ was conducted by Yangbo Ye \cite{Ye01}, but the necessary spectral interpretation was not available to complete the analysis.
With the uniform spectral expansion of \cite{HWI, HWII}, this type of analysis is now a feasible method to obtain the remaining Kloosterman zeta functions.
We leave these to future papers.

\section{Results}
\label{sect:Results}
As in \cite{BFG}, we define the exponential sum
\begin{align*}
	&S(m_1, m_2, n_1, n_2; D_1, D_2) \\
	& = \sum_{\substack{B_1, C_1 \summod{D_1}\\B_2, C_2 \summod{D_2}\\ }} \e{\frac{m_1B_1 + n_1(Y_1 D_2 - Z_1 B_2)}{D_1} + \frac{m_2B_2 + n_2(Y_2 D_1 - Z_2B_1)}{D_2}},
\end{align*}
where the sum is restricted to
\[ D_1C_2 + B_1B_2 + D_2C_1 \equiv 0 \summod{D_1D_2}, \qquad (B_1, C_1, D_1) = (B_2, C_2, D_2) = 1, \]
and the $Y_i$ and $Z_i$ are defined by
\[ Y_1B_1 + Z_1C_1 \equiv 1 \pmod{D_1}, \qquad Y_2B_2 + Z_2C_2 \equiv 1 \pmod{D_2}. \]
Then the long element $SL(3,\Z)$ Kloosterman sum is
\begin{align*}
	S_{w_l}(\psi_m,\psi_n;c) =& S(-n_2,-n_1,m_1,m_2;c_1, c_2).
\end{align*}
Note the difference in the definitions between \cite{SpectralKuz} and \cite{WeylI}, as explained in \cite{WeylI}*{Section 4}; we follow \cite{WeylI}.

\subsection{The arithmetic Kuznetsov formula}
Let $\mathcal{S}_2^d$ be bases of even or odd spherical $SL(2,\Z)$ Maass cusp forms as $d=0$ or $d=1$, respectively.
For $\phi \in \mathcal{S}_2^d$, denote its spectral parameters by $(\mu_\phi,-\mu_\phi)$ with $\mu_\phi \in i\R$ (by \cite{Roelcke01}) and the Hecke eigenvalues of the maximal parabolic Eisenstein series of weight $d$ attached to $\phi\in\mathcal{S}_2^d$ with the additional spectral parameter $r \in \C$ by $\lambda_\phi(m,r)$ (see \cite{HWI}*{Section 5.5}).

Let $\mathcal{S}_3^d$ be a basis of vector-valued $SL(3,\Z)$ Maass cusp forms of weight $d$.
For $\varphi\in\mathcal{S}_3^d$ we denote its Hecke eigenvalues by $\lambda_\varphi(m)$ and its spectral parameters by $\mu_\varphi$.
Denote the Hecke eigenvalues of the minimal parabolic spherical Eisenstein series by $\lambda_E(m,\mu)$ (see \cite{HWI}*{Section 4.2}).

The cuspidal and Eisenstein series parts of both $d=0,1$ (minimal-weight) spectral Kuznetsov formulae with $m,n\in\Z^2$, $m_1 m_2 n_1 n_2 \ne 0$ can be written
\begingroup\allowdisplaybreaks\begin{align}
\label{eq:KuzCuspPart}
	\mathcal{C}^d(F) =& \frac{2\pi}{3} \sum_{\varphi \in \mathcal{S}_3^d} F(\mu_\varphi) \frac{\wbar{\lambda_\varphi(m)} \lambda_\varphi(n)}{L(\Ad^2 \varphi,1)}, \\
\label{eq:MaxEisenPart}
	\mathcal{E}^d_\text{max}(F) =& \frac{1}{2 i} \sum_{\phi \in \mathcal{S}_2^d} \int_{\Re(r)=0} \frac{F(\mu_\phi+r,-\mu_\phi+r,-2r) \wbar{\lambda_\phi(m,r)} \lambda_\phi(n,r)}{L(\phi,1+3r) L(\phi,1-3r) L(\Ad^2 \phi,1)} dr, \\
\label{eq:MinEisenPart}
	\mathcal{E}^0_\text{min}(F) =& -\frac{1}{12 \pi} \int_{\Re(\mu)=0} \frac{F(\mu) \lambda_E(n, \mu) \wbar{\lambda_E(m, \mu)}}{\prod_{i<j} \abs{\zeta(1+\mu_i-\mu_j)}^2} d\mu,
\end{align}\endgroup
where $F$ is a smooth function of $\mu$ satisfying certain regularity and symmetry conditions; see Theorem \ref{thm:SpectralKuz}, and
\[ d\mu = d\mu_1 \, d\mu_2 = d\mu_1 \, d\mu_3 = d\mu_2 \, d\mu_3. \]

For a smooth, compactly supported function $f$ on $(\R^+)^2$, define the integral transforms
\begin{align}
\label{eq:F0def}
	F_0(\mu) = F_0(f;\mu) =& \frac{16}{\pi^4} \cosmu^0(\mu) \int_{(\R^+)^2} f(t) W^{0*}(t,-2\mu) \frac{dt_1 \, dt_2}{(t_1 t_2)^2}, \\
\label{eq:F1def}
	F_1(\mu) = F_1(f;\mu) =& \frac{1}{2} \tan\frac{\pi}{2}(\mu_1-\mu_3)\tan\frac{\pi}{2}(\mu_2-\mu_3) F_0(\mu),
\end{align}
where
\begin{align}
\label{eq:cos0}
	\cosmu^0(\mu) = \frac{2}{\pi} \prod_{i<j} \cos\frac{\pi}{2}(\mu_i-\mu_j),
\end{align}
and $W^{0*}(y,\mu)$ is the completed spherical Whittaker function
\begin{align}
\label{eq:SphWhitt}
	W^{0*}(y,\mu) =& \frac{1}{4\pi^2} \int_{-i\infty}^{i\infty} \int_{-i\infty}^{i\infty} (\pi y_1)^{1-s_1} (\pi y_2)^{1-s_2} G^0(s,\mu) \frac{ds_1 \, ds_2}{(2\pi i)^2}, \\
	G^0(s,\mu) =& \frac{\prod_{i=1}^3 \Gamma\paren{\frac{s_1-\mu_i}{2}} \Gamma\paren{\frac{s_2+\mu_i}{2}}}{\Gamma\paren{\frac{s_1+s_2}{2}}},
\end{align}
with spectral parameters $\mu\in\C^3$, $\mu_1+\mu_2+\mu_3=0$.
The contours of \eqref{eq:SphWhitt} should pass to the right of the poles of $G^0(s,\mu)$.
Notice that in case $\wbar{\mu}$ is a permutation of $-\mu$ (in particular if $\Re(\mu)=0$), we have $\wbar{W^{0*}(y,\mu)} = W^{0*}(y,-\mu)$.

The arithmetic Kuznetsov formula is given in terms of these integral transforms.
\begin{thm}
\label{thm:ArithKuz}
Suppose $m,n\in \Z^2$ such that $m_1 n_2, m_2 n_1 > 0$, and $f$ smooth and compactly supported on $(\R^+)^2$, then
\begin{align*}
	\mathcal{K}_L(f) :=& \sum_{c_1,c_2\in\N} \frac{S_{w_l}(\psi_m,\psi_n,c)}{c_1 c_2} f\paren{\tfrac{2\sqrt{m_1 n_2 c_2}}{c_1}, \tfrac{2\sqrt{m_2 n_1 c_1}}{c_2}}.
\end{align*}
has the spectral interpretation
\begin{align}
\label{eq:ArithKuzThm}
	\mathcal{K}_L(f) = \mathcal{C}^0(F_0)+\mathcal{C}^1(F_1)+\mathcal{E}^0_\text{max}(F_0)+\mathcal{E}^0_\text{min}(F_0)+2\mathcal{E}^1_\text{max}(F_1).
\end{align}
\end{thm}

Of course, the hypotheses on $f$ may be reduced to absolute convergence of $\mathcal{K}_L(f)$, absolute convergence of \eqref{eq:F0def}, and a bound of the form $F_0(\mu) \ll (1+\norm{\mu})^{-5-\epsilon}$, with $\norm{\mu}$ the usual Euclidean norm on $\C^3$, see Section \ref{sect:ArithKuzConv}.
An optimal result is well beyond the scope of this paper, but the following proposition gives a sufficient (though certainly not necessary) set of conditions:
\begin{prop}
\label{prop:ArithKuzConv}
Suppose $f:(\R^+)^2\to\C$ and that there exists $\sigma_1 > \theta$, $\sigma_2 \ge 0$, with $\theta$ as in Theorem \ref{thm:SmoothSums}, such that
\begin{enumerate}
\item $f(y)$ is bounded by $(y_1 y_2)^{1+\epsilon}$ as $y_1 y_2 \to 0$,
\item for $j_1,j_2\le 2N,j_1+j_2\le 2N$, the derivatives $\partial_{y_1}^{j_1} \partial_{y_2}^{j_2} f(y)$ exist and are continuous,
\item for $j_1,j_2\le 2N,j_1+j_2\le 2N-1$, the derivatives $y_1^{j_1} y_2^{j_2} \partial_{y_1}^{j_1} \partial_{y_2}^{j_2} (y_1 y_2 f(y))$ are bounded by $(y_1 y_2)^{1+\epsilon}$ as $y_1 y_2 \to 0$ and are otherwise polynomially bounded in the coordinates of $y$,
\item $(\wtilde{\Delta_1})^N (y_1 y_2 f(y)) \ll (y_1 y_2)^{1+\sigma_1} (y_1+y_2)^{2\sigma_2}$ for all $y\in(\R^+)^2$,
\end{enumerate}
where $N$ is the least integer strictly larger than $\sigma_1+\sigma_2+\frac{9}{4}$, and
\begin{align}
\label{eq:tildeDelta1Def}
	\wtilde{\Delta_1} = -y_1^2 \partial_{y_1}^2-y_2^2\partial_{y_2}^2+y_1 y_2 \partial_{y_1} \partial_{y_2}+4\pi^2(y_1^2+y_2^2).
\end{align}
Then \eqref{eq:ArithKuzThm} holds for the test function $f$.
\end{prop}

The strength of Theorem \ref{thm:SmoothSums} relies on non-trivial bounds on the spectral parameters; this was essentially proved in \cite{KS01}:
\begin{thm}[Kim, Sarnak]
\label{thm:KS}
	The spectral parameters $\mu$ of a cusp form of weight at most one satisfy $\abs{\Re(\mu_j)} < \frac{5}{14}$.
\end{thm}
Although the theorem given in \cite{LRS} is only stated for spherical forms, it trivially extends to the weight-one case; see the discussion in \cite[Section 9.1]{WeylI}.

Theorem \ref{thm:SmoothSums} now follows from Theorem \ref{thm:ArithKuz} and an analysis of the spectral expansion; the proof is given in Section \ref{sect:SmoothSums}, below.

\subsection{The Kloosterman zeta function}
\label{sect:KloosZeta}
At this point, one would like to give the meromorphic continuation of the unweighted Kloosterman zeta function.
We detour momentarily to analyze the failure of the obvious choice:
It seems that one should, with the usual limiting argument, be able to take $f(y)=(\pi y_1)^{2s_1} (\pi y_2)^{2s_2}$ in Theorem \ref{thm:ArithKuz} so that the Kloosterman sum side becomes
\[ \mathcal{K}_L(f) = (4 \pi^2 m_1 n_2)^{s_1} (4 \pi^2 m_2 n_1)^{s_2} \sum_{c_1,c_2\in\N} \frac{S_{w_l}(\psi_m,\psi_n,c)}{c_1^{1+2s_1-s_2} c_2^{1+2s_2-s_1}}, \]
which converges on
\[ \set{s\in\C^2 \setdiv \tfrac{1}{2} < \Re(s_1) < 2\Re(s_2) -\tfrac{1}{2}, \tfrac{1}{2} < \Re(s_2) < 2\Re(s_1) -\tfrac{1}{2}}. \]
If we do this, then $F_0(\mu)$ becomes
\begin{align*}
	F_0(\mu) =& \frac{4}{\pi^4} \cosmu^0(\mu) G^0(2s,-2\mu),
\end{align*}
and here is where our problems begin.
The function $F_0(\mu)$ has exponential decay unless one of the $\mu_i$ is small (see Section \ref{sect:Residues}), but if, say, $\mu_2$ is not too large (compared to $\Im(s)$), and $\mu_1 \asymp -\mu_3 \asymp T$ is large (which rules out the complementary spectrum), we see $F_0(\mu) \asymp T^{2\Re(s_1+s_2)-2}$ for $\Re(\mu)=0$ by Stirling's formula.
Unfortunately, the number of such cusp forms (arithmetically weighted) is within a constant multiple of $T^4$ (see \cite[Theorem 1]{Val01} for the spherical case), so the spectral expansion only converges absolutely for $2\Re(s_1+s_2) < -2$.
So it seems there is no overlap between the two representations, or more accurately, the usual limiting argument doesn't apply; in particular, it is not possible to choose $\sigma_1$ and $\sigma_2$ for Proposition \ref{prop:ArithKuzConv} due to the term $4\pi^2(y_1^2+y_2^2)$ in \eqref{eq:tildeDelta1Def}, which requires increasing $\sigma_1$ and $\sigma_2$ as the number of derivatives $N$ grows at precisely the same rate.
This differs from the behavior of the spectral expansions of the Kloosterman zeta functions on $GL(2)$ which have exponential convergence in the spectral expansion \cite{Moto01,Kiral01}, and we note that on the generic region where all $\abs{\mu_i} \asymp T$, this naive test function would also produce exponential convergence.

The solution on $GL(2)$ is to use an exponential to control the behavior at infinity, and write the weighted Kloosterman zeta function as a sum of shifts of the unweighted Kloosterman zeta function \cite{Sel04}.
The meromorphic continuation of the unweighted function can then be obtained from that of the weighted function in stages.
This iterative approach was formalized by Kiral \cite{Kiral01}, who obtained a sort of M\"obius inversion for the power series expansion of the exponential to express the unweighted Kloosterman zeta function as a sum of shifts of the weighted Kloosterman zeta function, in the opposite-sign case.

To that end, one might choose
\[ f(y) = (\pi y_1)^{2s_1} (\pi y_2)^{2s_2} \exp\paren{-(\pi^9 y_1^5 y_2^4)^2-(\pi^9 y_1^4 y_2^5)^2}. \]
Then Proposition \ref{prop:ArithKuzConv} applies with
\[ \sigma_1=\theta+\epsilon, \qquad \sigma_2=\epsilon, \qquad N=3 \]
(not an optimal choice) on the region
\[ 10 s_1-8s_2,10 s_2-8s_1 > 24+\theta. \]
Unfortunately, it can be seen that the spectral side for this new weight function does not converge for, say $-\frac{1}{2} < \Re(s_1) < 0$ and $\Re(s_2) > -\frac{1}{2}$, for exactly the same reason as the naive test function in the previous attempt.
In fact, studying the failure of this example leads one to believe that any test function of the form
\[ f(y) = (\pi y_1)^{2s_1} (\pi y_2)^{2s_2} f_0(y) \]
where $f_0(y)$ satisfies some conditions on its Mellin transform (e.g. existence) and
\[ f_0(y) = 1 + O\paren{\sum_i y_1^{\alpha_i} y_2^{\beta_i}} \]
as $y \to 0$ with $\alpha_i$ and $\beta_i$ strictly positive, will encounter a similar difficulty.
This is unfortunate, as we need $2\alpha_i-\beta_i \ge 0, 2\beta_i-\alpha_i \ge 0$ with one inequality strict -- that is, the exponents of both $c_1$ and $c_2$ must be non-positive with one exponent strictly negative -- in each error term to apply the iterative method.

Thus we consider the simple weighted Kloosterman zeta function of the theorem.
We collect the relevant analysis into the following proposition:
\begin{prop}
\label{prop:KloosZetaTechnical}
For $f_s$, $\tilde{s}$ and $Z^*_{m,n}(s)$ as in Theorem \ref{thm:KloosZeta},
\begin{enumerate}
\item \eqref{eq:ArithKuzThm} holds with $\mathcal{K}_L(f_s) = Z^*_{m,n}(s)$ on the region $2\Re(\tilde{s}_1),2\Re(\tilde{s}_2) > 1$,

\item $F_0(s,\mu) := F_0(f_s,\mu)$ and $F_1(s,\mu):=F_1(f_s,\mu)$ extend to meromorphic functions of all $s$ and $\mu$,

\item when $\abs{\tilde{s}} < T$, $\Re(\mu)$ is in some compact set $M$, and
	\[ \abs{\tilde{s}_1+\mu_j+\ell},\abs{\tilde{s}_2-\mu_j+\ell} > \eta > 0, \qquad \text{ for each $j$ and all }\ell\in\N_0, \]
	we have $F_0(s,\mu),F_1(s,\mu) \ll_{T,M,\eta} (1+\norm{\mu})^{-100}$,

\item when the coordinates of $\mu$ are distinct modulo $\Z$, $F_0(s,\mu)$ and $F_1(s,\mu)$ have simple poles whenever $\tilde{s}_1=-\mu_j-\ell$ or $\tilde{s}_2=\mu_j-\ell$, $\ell\in\N_0$,

\item when $\mu=(it,it,-2it)$ with $0 \ne t\in\R$, $F_0(s,\mu)$ and $F_1(s,\mu)$ have double poles whenever $\tilde{s}_1=-it-\ell$ or $\tilde{s}_2=it-\ell$, and simple poles whenever $\tilde{s}_1=2it-\ell$ or $\tilde{s}_2=-2it-\ell$, $\ell\in\N_0$,

\item $F_0(s,0)$ has triple poles whenever $\tilde{s}_1=-\ell$ or $\tilde{s}_2=-\ell$, $\ell\in\N_0$,

\item $F_1(s,0)=0$.
\end{enumerate}
\end{prop}
These computations are of similar type to \cite{IwSMAF}*{Section 9.2} and \cite{BFG}, and we put them off to Sections \ref{sect:Residues} and \ref{sect:KloosZetaTechnical}.
Notice the order of each pole is the maximum number of terms of \eqref{eq:KloosZetaPoles} which coincide for any given $\mu$.

From parts 1,2 and 3 of the proposition, the right-hand side of \eqref{eq:ArithKuzThm} gives a meromorphic continuation of $Z^*_{m,n}(s)$ to all $s$.
The spectral parameters $\mu$ of spherical and weight-one forms must satisfy the unitaricity condition and the bound towards Ramanujan-Selberg (e.g. $\theta=\frac{5}{14}$ as in Theorem \ref{thm:KS}):
\begin{align}
\label{eq:muProps}
	\textstyle \max_i \abs{\Re(\mu_i)} < \theta<\frac{1}{2}, \qquad -\wbar{\mu} \text{ a permutation of } \mu.
\end{align}
For such $\mu$, having $\mu_j-\mu_k\in\Z$ implies $\mu$ is a permutation of $(it,it,-2it)$, but $F_0(s,\mu)$ is permutation-invariant in $\mu$ and $F_1(s,\mu)$ is zero on the other permutations, so the proposition covers all relevant cases in $\mu$.
Thus Theorem \ref{thm:KloosZeta} follows immediately from Proposition \ref{prop:KloosZetaTechnical} and the associated Weyl laws \cite{Val01}*{Theorem 1} and \cite{WeylI}*{Theorem 1}.

It is known \cite{CHJT,MillerThesis,Miller02} that $\mu=0$ does not correspond to a spherical cusp form on $SL(3,\Z)$, so we only need to worry about $\mu=0$ in the terms $\mathcal{E}^0_\text{max}$ and $\mathcal{E}^0_\text{min}$.
On the other hand, it is known that $0$ is not the spectral parameter of any cusp form on $SL(2,\Z)$ \cite{Roelcke01}, i.e. in \eqref{eq:MaxEisenPart} we have $\mu_\phi \ne 0$, so the $\mu$ argument of $F_0(s,\mu)$ is never zero in $\mathcal{E}^0_\text{max}$, and the integrand of $\mathcal{E}^0_\text{min}$ has a triple zero at $\mu=0$ due to the poles of the six zeta functions in the denominator there.
So the point $\mu=0$ and its triple pole need not enter into our computations.

The residues of $Z_{m,n}^*(s)$ at say $\tilde{s}_1=-\mu_k-\ell$ for $\mu$ the spectral parameters of some cusp form (of either weight) and $\ell\in\N_0$ are then
\begin{align}
	\frac{2\pi}{3} \sum_{d\in\set{0,1}} \sum_{j\in\set{1,2,3}} \sum_{\substack{\varphi \in \mathcal{S}_3^d\\ \mu_{\varphi,j}=\mu_k}} \frac{\wbar{\lambda_\varphi(m)} \lambda_\varphi(n)}{L(\Ad^2 \varphi,1)} \res_{\tilde{s}_1=-\mu_k-\ell} F_d(s, \mu_\varphi)+\text{Eisenstein contribution}.
\end{align}
The residues of $F_0(s,\mu)$ (from which follow those of $F_1(s,\mu)$), up to symmetry and away from $\mu=0$, are listed in \eqref{eq:F0fromTildeF0},\eqref{eq:tildeF0ResGenMu},\eqref{eq:tildeF0ResDegenMu}.

We have omitted the polar analysis of the minimal and maximal parabolic Eisenstein terms, as they are somewhat more complicated.
Essentially, when $\tilde{s}$ is close to a pole of $F_d(s,\mu)$, we deform the $\mu$ or $r$ contour past the pole, then ask if $\tilde{s}$ is a pole of the resulting residual term.
The essential idea is the same as writing
\begin{align*}
	y^s e^{-y} =& \begin{cases}\displaystyle \int_{\Re(u)=\epsilon} \Gamma(u+s) y^{-u} \frac{du}{2\pi i} & \Re(s) \ge 0, \\
	\displaystyle \int_{\Re(u)=\epsilon} \Gamma(u+s) y^{-u} \frac{du}{2\pi i}+\sum_{\ell=0}^T \frac{(-1)^\ell}{\ell!} y^{s+\ell} & -T-1\le\Re(s) < -T \in -\N_0.\end{cases}
\end{align*}
The summands of the $\ell$ sum here are referred to as the residual terms.
(Note that the $\epsilon$ in the contour necessarily depends on $\Re(s)$.)

Poles of the residual terms may result from the zeros of the $L$-functions in the denominators of \eqref{eq:MaxEisenPart} or \eqref{eq:MinEisenPart}, or from secondary poles of the residue of $F_d(s,\mu)$.
This is further complicated by the fact that it is not possible to manipulate the $r$ or $\mu$ contours so that all three coordinates of $\mu$ lie on the same side of $\Re(\mu)=0$ (since the coordinates sum to zero), so these secondary poles have a different form when one or more of $\tilde{s}_i \in -\N_0+i\R$.

The possible secondary poles of the residual terms coming from the minimal parabolic Eisenstein series occur whenever $\Re(\tilde{s}_1),\Re(\tilde{s}_2) \le 0$ and $\tilde{s}_1+\tilde{s}_2 \in -\N_0$.
The possible poles coming from the zeta functions in the denominator occur whenever $\Re(\tilde{s}_1) \le -\ell_1 \in -\N_0$, $\Re(\tilde{s}_2) \le -\ell_2 \in -\N_0$, and any of the numbers
\begin{align}
\label{eq:PoleNumbers}
	\tilde{s}_1+\ell_1+\tilde{s}_2+\ell_2, \qquad 2(\tilde{s}_1+\ell_1)-(\tilde{s}_2+\ell_2), \qquad 2(\tilde{s}_2+\ell_2)-(\tilde{s}_1+\ell_1)
\end{align}
are of the form $\pm(1-z)$ for $z$ any non-trivial zero of $\zeta(z)$ (the trivial zeros are killed by $\cosmu^0(\mu)$).

The possible secondary poles of the residual terms coming from the maximal parabolic Eisenstein series attached to $\phi\in\mathcal{S}^d_2$ occur whenever $\Re(\tilde{s}_1) \le -\ell_1 \in -\N_0$, $\Re(\tilde{s}_2) \le -\ell_2 \in -\N_0$ and either $\tilde{s}_1+\tilde{s}_2 \in -\N_0$ or any of the numbers \eqref{eq:PoleNumbers} is $\pm 2\mu_\phi$.
The possible poles coming from the $L$-functions in the denominator occur whenever one of $\Re(\tilde{s}_i) \le -\ell_i \in -\N_0$ and either of the numbers 
\[ \tfrac{1}{2}(\tilde{s}_i+\ell_i), \qquad \tilde{s}_i+\ell_i\pm \mu_\phi, \]
are of the form $\pm\frac{1-z}{3}$ for $z$ any zero of $L(\phi,z)$.

\section{Acknowledgements}
The author would like to thank Valentin Blomer for many helpful comments.
The author would also like to thank the referees for their comments which uncovered a significant mistake with the Kloosterman zeta function in a previous version of this paper.

\section{The spherical Whittaker function}
We follow the notation of \cite{HWI}*{Section 2.1}.
In particular, let $U(\R)$ be the space of upper-triangular, unipotent matrices in $PSL(3,\R)$ and $Y^+$ the space of positive diagonal matrices, in the form $y=\diag(y_1 y_2, y_1, 1)$ with measure $dy=\frac{dy_1 \, dy_2}{(y_1 y_2)^3}$.
We frequently treat $y\in Y^+$ as elements of $(\R^+)^2$ as the multiplication is the same, and the same for the signed diagonal matrices $Y \cong (\R^\times)^2$.

If we take the usual power function $p_{\rho+\mu}(y) = \abs{y_1}^{1-\mu_3} \abs{y_2}^{1+\mu_1}$ with $\rho=(1,0,-1)$, then the incomplete spherical Whittaker function is defined by the Jacquet integral
\[ W^0(g,\mu,\psi_n) = \int_{U(\R)} p_{\rho+\mu}(w_l u g) \wbar{\psi_n(u)} du, \]
for $\psi_n$ a character of $U(\R)$.
It can be shown (see \cite{Gold01} or the more complicated cases of \cite{HWII}) that this has the Mellin-Barnes integral \eqref{eq:SphWhitt} with
\begin{gather*}
	\Lambda(\mu) = \pi^{-\frac{3}{2}+\mu_3-\mu_1} \Gamma\paren{\tfrac{1+\mu_1-\mu_2}{2}} \Gamma\paren{\tfrac{1+\mu_1-\mu_3}{2}} \Gamma\paren{\tfrac{1+\mu_2-\mu_3}{2}}, \\
	W^{0*}(y,\mu) := \Lambda(\mu)W^0(y,\mu,\psi_{1,1}).
\end{gather*}
This fixes the various normalization constants to give us a consistent starting point.

The power function induces an action of the Weyl group $W$ (as in \cite{HWI}*{Section 2.1}) on the coordinates of $\mu$ by $p_{\mu^w}(y)=p_\mu(wyw^{-1})$ for $w \in W$:
\begin{align*}
	\mu^I =& \paren{\mu_1,\mu_2,\mu_3}, & \mu^{w_2} =& \paren{\mu_2,\mu_1,\mu_3}, & \mu^{w_3} =& \paren{\mu_1,\mu_3,\mu_2}, \\
	\mu^{w_4} =& \paren{\mu_3,\mu_1,\mu_2}, & \mu^{w_5} =& \paren{\mu_2,\mu_3,\mu_1}, & \mu^{w_l} =& \paren{\mu_3,\mu_2,\mu_1}.
\end{align*}
We say that a function $F(\mu)$ is Weyl invariant if $F(\mu)=F(\mu^w)$ for all $w\in W$.

It is known that the Whittaker functions give a continuous basis of the Schwartz-class functions on $(\R^+)^2$ by the Kontorovich-Lebedev theorem of Wallach \cite{Wallach} (see \cite{GoldKont}):
\begin{thm}
\label{thm:KontLeb}
Let $F(\mu)$ be Schwartz-class, Weyl-invariant and holomorphic for $\Re(\mu)$ in a neighborhood of $0$, and let $f(y)$ be Schwartz class on $Y^+$, then
\begin{align*}
	F(\mu) =& \int_{Y^+} \int_{\Re(\mu')=0} F(\mu') W^{0*}(y,\mu') \sinmu^0(\mu') d\mu'\,W^{0*}(y,-\mu) dy, \\
	f(y) =& \int_{\Re(\mu)=0} \int_{Y^+} f(t) W^{0*}(t,-\mu) dt \, W^{0*}(y,\mu) \sinmu^0(\mu) d\mu, 
\end{align*}
with the integrals to be taken iteratively, where
\begin{align}
\label{eq:sin0}
	\sinmu^0(\mu) := \frac{1}{6 (2\pi i)^2 \prod_{i\ne j} \Gamma\paren{\tfrac{\mu_i-\mu_j}{2}}} = \frac{1}{192\pi^5} \prod_{i<j} (\mu_i-\mu_j)\sin\frac{\pi}{2}(\mu_i-\mu_j).
\end{align}
\end{thm}
Note: In comparing to \cite{GoldKont}*{Definitions 1.1, 1.2 and Theorem 1.3}, we have the conversion $3it=(\mu_1-\mu_2,\mu_2-\mu_3)$ and Stade's formula \cite{GoldKont}*{eq. (2.1)} in that paper is missing a factor $2^{n-2}$ from the denominator of the right-hand side.
The conjugations may be explained by $\wbar{W^{0*}(y,\mu)}=W^{0*}(y,-\mu)$ on $\Re(\mu)=0$.
The method of proof in \cite{GoldKont} requires holomorphy of $F$, and they technically only prove one direction, but Wallach's proof culminating in \cite{Wallach}*{Section 15.10.3 eqs. (1) and (2)} is somewhat daunting.
The constants and the spectral measure are verified through the computation of Section \ref{sect:SphericalConstants}.

\section{Spectral Kuznetsov Formulae}
We refer to \cite{WeylI}*{Section 4.5} for the definitions of the power-series solutions $J_w(y,\mu)$ to the differential equations satisfied by the Kuznetsov kernel functions.
Their precise definition is not needed here, as we will only use them in an algebraic sense.

For $F(\mu)$ Schwartz-class and holomorphic for $\Re(\mu)$ in a neighborhood of $0$, we define the integral transformation
\begin{align*}
	H_w^d(F; y) =& \frac{1}{\abs{y_1 y_2}} \int_{\Re(\mu)=0} F(\mu) K^d_w(y, \mu) \specmu^d(\mu) d\mu
\end{align*}
using the spectral weights
\begin{align}
\label{eq:spec0}
	\specmu^0(\mu) :=& \frac{\sinmu^0(\mu)}{\cosmu^0(\mu)} = \frac{1}{384\pi^4} \prod_{i<j} (\mu_i-\mu_j)\tan\frac{\pi}{2}(\mu_i-\mu_j), \\
\label{eq:spec1}
	\specmu^1(\mu) =& \frac{1}{64\pi^4} \paren{\prod_{i<j} (\mu_i-\mu_j)} \cot\frac{\pi}{2}(\mu_1-\mu_3) \cot\frac{\pi}{2}(\mu_2-\mu_3) \tan\frac{\pi}{2}(\mu_1-\mu_2),
\end{align}
and the kernel functions
\begingroup\allowdisplaybreaks\begin{align}
\label{eq:KwIdDef}
	K_I^d(y,\mu) =& 1, \\
\label{eq:Kw40Def}
	K_{w_4}^0(y,\mu) =& \frac{1}{8\pi} \sum_{w\in W_3} \frac{J_{w_4}(y,\mu^w)}{\sin \frac{\pi}{2}(\mu^w_1-\mu^w_3) \sin \frac{\pi}{2}(\mu^w_2-\mu^w_3)}, \\
\label{eq:Jw41Def}
	J^1_{w_4}(y,\mu) =& -\sin\frac{\pi}{2}(\mu_1-\mu_2) J_{w_4}(y,\mu) - i\varepsilon_1 \cos\frac{\pi}{2}(\mu_1-\mu_3) J_{w_4}(y,\mu^{w_4}), \\
	& \qquad + i\varepsilon_1 \cos\frac{\pi}{2}(\mu_2-\mu_3) J_{w_4}(y,\mu^{w_5}), \nonumber \\
\label{eq:Kw41Def}
	K^1_{w_4}(y,\mu) =& \frac{1}{8\pi} \frac{J^1_{w_4}(y,\mu)}{\cos\frac{\pi}{2}(\mu_1-\mu_3)\cos\frac{\pi}{2}(\mu_2-\mu_3)\sin\frac{\pi}{2}(\mu_1-\mu_2)}, \\
\label{eq:Kw5dDef}
	K_{w_5}^d(y,\mu) =& K_{w_4}^d((-y_2,y_1),-\mu), \\
\label{eq:Kwl0Def}
	K_{w_l}^0(y,\mu) =& \frac{1}{16 \pi} \frac{\sum_{w\in W_3} \paren{J_{w_l}(y,\mu^{w w_2})-J_{w_l}(y,\mu^w)}}{\prod_{i<j} \sin \frac{\pi}{2}\paren{\mu_i-\mu_j}}, \\
	J^1_{w_l}(y,\mu) =& \varepsilon_2 J_{w_l}(y,\mu)+\varepsilon_1 J_{w_l}(y,\mu^{w_4})+\varepsilon_1 \varepsilon_2 J_{w_l}(y,\mu^{w_5}), \\
\label{eq:Kwl1Def}
	K^1_{w_l}(y,\mu) =& -\frac{1}{16\pi} \frac{J^1_{w_l}(y,\mu)-J^1_{w_l}(y,\mu^{w_2})}{\cos\frac{\pi}{2}(\mu_1-\mu_3)\cos\frac{\pi}{2}(\mu_2-\mu_3)\sin\frac{\pi}{2}(\mu_1-\mu_2)},
\end{align}\endgroup
with $\varepsilon=\sgn(y)$, and $W_3=\set{I,w_4,w_5}$ the subgroup of order 3 permutations in the Weyl group.
The weight-one formulae are \cite{WeylI}*{(4)-(7)} and the surrounding displays; the spherical formulae are derived in Section \ref{sect:SphericalConstants}, below.

With these integral transforms, we consider the following averages of Kloosterman sums at each Weyl cell with $m,n\in\Z^2$, $m_1 m_2 n_1 n_2 \ne 0$:
\begingroup\allowdisplaybreaks\begin{align*}
	\mathcal{K}_I^d(F) =& \delta_{\substack{\abs{m_1}=\abs{n_1}\\ \abs{m_2}=\abs{n_2}}} H_I^d(F;(1, 1)) \nonumber \\
	\mathcal{K}_4^d(F) =& \sum_{\varepsilon\in\set{\pm1}^2} \sum_{\substack{c_1,c_2\in\N\\ \varepsilon_1 m_2 c_1=n_1 c_2^2}} \frac{S_{w_4}(\psi_m,\psi_{\varepsilon n},c)}{c_1 c_2} H_{w_4}^d\paren{F; \paren{\varepsilon_1 \varepsilon_2 \tfrac{m_1 m_2^2n_2}{c_2^3n_1},1}}, \\
	\mathcal{K}_5^d(F) =& \sum_{\varepsilon\in\set{\pm1}^2} \sum_{\substack{c_1,c_2\in\N\\ \varepsilon_2 m_1 c_2= n_2 c_1^2}} \frac{S_{w_5}(\psi_m,\psi_{\varepsilon n},c)}{c_1 c_2} H_{w_5}^d\paren{F; \paren{1,\varepsilon_1 \varepsilon_2 \tfrac{m_1^2 m_2 n_1}{c_1^3 n_2}}}, \\
	\mathcal{K}_L^d(F) =& \sum_{\varepsilon\in\set{\pm1}^2} \sum_{c_1,c_2\in\N} \frac{S_{w_l}(\psi_m,\psi_{\varepsilon n},c)}{c_1 c_2} H_{w_l}^d\paren{F; \paren{\varepsilon_2 \tfrac{m_1 n_2 c_2}{c_1^2}, \varepsilon_1 \tfrac{m_2 n_1 c_1}{c_2^2}}},
\end{align*}\endgroup
and at last we are ready for the spectral Kuznetsov formulae:
\begin{thm}
\label{thm:SpectralKuz}
	Let $F(\mu)$ be Schwartz-class and holomorphic on $\set{\mu|\abs{\Re(\mu_i)} < \frac{1}{2}+\delta}$ for some $\delta > 0$.
	Suppose $m,n\in\Z^2$ with $m_1 m_2 n_1 n_2 \ne 0$ and recall \eqref{eq:KuzCuspPart}-\eqref{eq:MinEisenPart}.
	\begin{enumerate}
	\item[0.] If $F(\mu)$ is invariant under all $\mu \mapsto \mu^w$ and $F(\mu)=0$ whenever $\mu_i-\mu_j=\pm1, i\ne j$, then we have the spherical Kuznetsov formula,
	\[ \mathcal{C}^0(F)+\mathcal{E}^0_\text{max}(F)+\mathcal{E}^0_\text{min}(F)=\mathcal{K}_I^0(F)+\mathcal{K}_4^0(F)+\mathcal{K}_5^0(F)+\mathcal{K}_L^0(F). \]
	\item[1.] If $F(\mu)$ is invariant under the transposition $\mu_1 \leftrightarrow \mu_2$ and $F(\mu)=0$ whenever $\mu_1-\mu_2=\pm1$, then we have the weight-one Kuznetsov formula,
	\[ \mathcal{C}^1(F)+2\mathcal{E}^1_\text{max}(F) =\mathcal{K}_I^1(F)+\mathcal{K}_4^1(F)+\mathcal{K}_5^1(F)+\mathcal{K}_L^1(F). \]
	\end{enumerate}
\end{thm}
In the weight-one case, this is \cite{WeylI}*{Theorem 6}; the spherical case is \cite{SpectralKuz}*{Theorems 1 and 4}, with the corrected constants which we now compute.

\subsection{The constants in the spherical formula}
\label{sect:SphericalConstants}
In what follows, the leading constants in the terms of the arithmetic and spectral sides of both spectral Kuznetsov formulae will be essential to the construction of the arithmetic Kuznetsov formula.
Unfortunately, these constants in \cite{SpectralKuz} are generally incorrect.
The key difficulties are that the paper \cite{SpectralKuz} relied on outside formulae such as spherical inversion (section 1.8 there; the constant there is extremely difficult to verify), the Mellin transform of the Whittaker function (section 1.4 there), Li's Kuznetsov formula (section 1.9 there) as computed in \cite{Gold01} and \cite{MeThesis}, especially Stade's formula (eq. (17) there), that \cite{SpectralKuz} has the extra step of computing the Fourier transform of the spherical function (eq. (22) there), and that \cite{SpectralKuz} attempts to evaluate the final $x'$ and $t$ integrals simultaneously (see the comment at the bottom of page 5 in \cite{WeylI}).

The papers \cite{WeylI} and \cite{WeylII} avoid these difficulties by computing all of the formulae directly, and there are far fewer steps.
Most of these computational errors in the constants of \cite{SpectralKuz} were eventually corrected, to the point that the ratio of the constants on the trival term of the arithmetic side and cuspidal term of the spectral side is known to be correct by heuristic verification (see \cite{SatoTateLaw}), but the remaining terms have no such means of verification, and it is difficult to check that the various corrections were properly carried through to those terms.
The solution is to repeat the computation from first principles using the methods of \cite{WeylI} and \cite{WeylII}.
We briefly summarize the computation here as it is essentially identical to those papers.

Stade's formula is computed (from \eqref{eq:SphWhitt} as in \cite{WeylI}*{Section 5}) to be
\begin{align*}
	\Psi^0(\mu,\mu',t) :=& \int_{Y^+} W^{0*}(y,\mu) W^{0*}(y,\mu') (y_1^2 y_2)^t dy = \frac{1}{4\pi^{3t} \Gamma\paren{\tfrac{3t}{2}}} \prod_{i,j} \Gamma\paren{\tfrac{t+\mu_i+\mu_j'}{2}},
\end{align*}
and at $t=1$ and $t=0$, this is
\[ \frac{1}{\cosmu^0(\mu)} := \Psi^0(\mu,-\mu,1), \qquad \frac{1}{\sinmu^0(\mu)} = 2(2\pi i)^2 \lim_{t\to 0} t^2 \Psi^0(\mu,-\mu,t), \]
which produce the formulae \eqref{eq:cos0}, \eqref{eq:sin0} and \eqref{eq:spec0}.

The generalization of Kontorovich-Lebedev inversion was already stated in Theorem \ref{thm:KontLeb}.

As in \cite{WeylII}*{Section 7.2} (see also the last paragraph of the introduction to \cite{WeylI}*{Section 7.1}), the kernel functions are defined by the Riemann integral
\begin{align*}
	K_w^0(y,\mu) W^{0*}(t,\mu) =& \int_{\wbar{U}_w(\R)} W^{0*}(ywxt,\mu) \,\wbar{\psi_{1,1}(x)} dx, \qquad y\in Y, t\in Y^+,
\end{align*}
and we determine the combinations of the power-series solutions by the asymptotics (via \eqref{eq:SphWhitt}) 
\begin{align*}
	W^{0*}(y,\mu) \sim& \frac{1}{2\pi} \sum_{w\in W_3} \pi^{-\mu^w_3} \abs{y_1}^{1-\mu^w_3} \Gamma\paren{\tfrac{\mu^w_3-\mu^w_1}{2}} \Gamma\paren{\tfrac{\mu^w_3-\mu^w_2}{2}} \int_{-i\infty}^{i\infty} (\pi y_2)^{1-s_2} \Gamma\paren{\tfrac{s_2+\mu^w_1}{2}} \Gamma\paren{\tfrac{s_2+\mu^w_2}{2}} \frac{ds_2}{2\pi i} \\
	\sim& \sum_{w\in W} \pi^{\mu^w_1-\mu^w_3} p_{\rho+\mu^w}(y) \Gamma\paren{\tfrac{\mu^w_3-\mu^w_1}{2}} \Gamma\paren{\tfrac{\mu^w_3-\mu^w_2}{2}} \Gamma\paren{\tfrac{\mu^w_2-\mu^w_1}{2}}
\end{align*}
as $y_1 \to 0$ in the first case or $y \to 0$ in general in the second, for $\mu$ with $\Re(\mu)=0$ and $\mu_i\ne\mu_j$, $i\ne j$.
For the long-element term inserting the $y\to 0$ asymptotic into the definition gives
\begin{align*}
	K_{w_l}^0(y,\mu) W^{0*}(t,\mu) \sim& \sum_{w\in W} \pi^{\mu^w_1-\mu^w_3} p_{\rho+\mu^w}(y) \frac{\Gamma\paren{\tfrac{\mu^w_3-\mu^w_1}{2}} \Gamma\paren{\tfrac{\mu^w_3-\mu^w_2}{2}} \Gamma\paren{\tfrac{\mu^w_2-\mu^w_1}{2}}}{\Lambda(\mu^w)} W^{0*}(t,\mu),
\end{align*}
and inserting the $y_1\to 0$ asymptotic gives
\begin{align*}
	K_{w_4}^0(y,\mu) W^{0*}(t,\mu) \sim& \sum_{w\in W_3} \pi^{1-3\mu^w_3} \abs{y_1}^{1-\mu^w_3} \frac{\Gamma\paren{\tfrac{\mu^w_3-\mu^w_1}{2}} \Gamma\paren{\tfrac{\mu^w_3-\mu^w_2}{2}}}{\Gamma\paren{\tfrac{1+\mu^w_2-\mu^w_3}{2}} \Gamma\paren{\tfrac{1+\mu^w_1-\mu^w_3}{2}}} W^{0*}(t,\mu),
\end{align*}
along $y_2=1$.
Then comparing with the first-term asymptotics of each $J_w(y,\mu)$ determines the kernel functions completely, giving \eqref{eq:Kw40Def} and \eqref{eq:Kwl0Def}.

The formulas \eqref{eq:KwIdDef} and \eqref{eq:Kw5dDef} for $K_I^0(y,\mu)$ and $K_{w_5}^0(y,\mu)$ follow from $\wbar{U}_I(\R)=\set{I}$ and $W^{0*}(y,\mu)=W^{0*}((y_2,y_1),-\mu)$ with the reasoning of \cite{WeylI}*{Section 7.1.4}.
This completes the corrections for the constants on the transforms $H_w$ in \cite{SpectralKuz}*{Theorem 1} as well as the formulae for the functions $K_w$ on page 6686 of \cite{SpectralKuz}, keeping in mind the notational differences discussed in the introduction to \cite{WeylI}*{Section 4}.

The constants on the spectral side of the spherical formula come from the correction
\begin{align*}
	\int_{\R^2} \Im(z_2)^{\frac{1}{2} \mp\mu_\phi} \paren{u_3^2+u_2^2+1}^{-\frac{3}{4}-\frac{3}{2}\mu_1} \, du_2 \, du_3 &= B\paren{\tfrac{1}{2},\tfrac{3\mu_1\mp\mu_\phi}{2}} B\paren{\tfrac{1}{2},\tfrac{3\mu_1\pm\mu_\phi}{2}},
\end{align*}
which brings \cite{SpectralKuz}*{Appendix A} in line with the computations of \cite{HWI}*{Section 5.5}.
In terms of $\cosmu^0(\mu_\varphi)$, the conversion to Hecke eigenvalues given in \cite{SpectralKuz}*{Appendix A} becomes
\begin{gather*}
	\cosmu^0(\mu_\varphi) \abs{\rho_\varphi^*(1)}^{-2} = \frac{3}{2\pi} L(\Ad^2 \varphi,1), \\
	\cosmu^0(\mu) \abs{\rho(1, \mu)}^{-2} = \frac{1}{8\pi} \abs{\zeta(1+\mu_1-\mu_2)\zeta(1+\mu_2-\mu_3)\zeta(1+\mu_1-\mu_3)}^2, \qquad \Re(\mu)=0, \\
	\cosmu^0(\mu_\phi+r,-\mu_\phi+r,-2r) \abs{\rho_\phi(1, r)}^{-2} = \frac{1}{\pi} L(\Ad^2 \phi,1) \abs{L(\phi,1+3r)}^2, \qquad \Re(r)=0.
\end{gather*}

\section{The Arithmetic Kuznetsov Formula}
We rely fundamentally on the fact \cite{SpectralKuz}*{Theorem 2} that
\begin{align}
\label{eq:KwlWhitt}
	K^0_{w_l}(y, \mu) =& \frac{\pi^5}{2} \cosmu^0(\mu) \sqrt{y_1 y_2} W^{0*}((2\sqrt{y_1},2\sqrt{y_2}),2\mu),
\end{align}
when $\sgn(y)=(1,1)$.

Suppose $F_0(\mu)$ is given by \eqref{eq:F0def} for some smooth, compactly supported $f:(\R^+)^2\to\C$, and define $F_1(\mu)$ by \eqref{eq:F1def}.
Notice that $F_0(\mu)$ satisfies the conditions of the spherical spectral Kuznetsov formula, Theorem \ref{thm:SpectralKuz}.0 since $W^{0*}(t,-2\mu)$ is entire in $\mu$ and an eigenfunction of the restricted Laplacian \eqref{eq:tildeDelta1Def} (using the positive Laplacian as in \cite{HWII}*{eqs. (42),(43)}):
\begin{align}
\label{eq:W0eigenfunc}
	\lambda_1(\mu) W^{0*}(y,\mu) &= \left.\Delta_1 W^{0*}(xyk,\mu)\right|_{x=k=I} = \wtilde{\Delta_1} W^{0*}(y,\mu), & \lambda_1(\mu) =& 1-\tfrac{\mu_1^2+\mu_2^2+\mu_3^2}{2}.
\end{align}
Then $F_1(\mu)$ also satisfies the conditions of the weight-one spectral Kuznetsov formula, Theorem \ref{thm:SpectralKuz}.1.
The zeros hypotheses are covered by the $\cosmu^0(\mu)$ factor and the symmetries by the Weyl-invariance of $W^{0*}(t,-2\mu)$.

Note that, in general, for functions $F_0$ and $F_1$ related by \eqref{eq:F1def}, even if $F_0$ is positive on the spherical spectrum, $F_1$ is not necessarily positive on the weight-one spectrum, and visa versa.

We add the two Kuznetsov formulae.

\subsection{The long-element term}
With $\varepsilon=\sgn(y)$, the long element weight functions can be written
\begin{align*}
	H_{w_l}^0(F_0; y) =& -\frac{6}{\abs{y_1 y_2}} \int_{\Re(\mu)=0} F_0(\mu) J_{w_l}(y, \mu) \frac{\specmu^0(\mu)}{16\pi\prod_{i<j} \sin \frac{\pi}{2}\paren{\mu_i-\mu_j}} d\mu, \\
	H_{w_l}^1(F_1; y) =& -\frac{6}{\abs{y_1 y_2}} \int_{\Re(\mu)=0} \paren{\varepsilon_2 +\varepsilon_1 +\varepsilon_1 \varepsilon_2} F_0(\mu) J_{w_l}(y, \mu) \frac{\specmu^0(\mu)}{16\pi\prod_{i<j} \sin \frac{\pi}{2}\paren{\mu_i-\mu_j}} d\mu,
\end{align*}
which follows directly from \eqref{eq:Kwl0Def} and \eqref{eq:Kwl1Def} using $F_0(\mu^w)=F_0(\mu)$ for all $w\in W$.
Combining these expresions into $H_{w_l}^*(F_0; y) := H_{w_l}^0(F_0; y)+H_{w_l}^1(F_1; y)$ with $1+\varepsilon_2 +\varepsilon_1 +\varepsilon_1 \varepsilon_2=4\delta_{\varepsilon=(1,1)}$ gives
\begin{align*}
	H_{w_l}^*(F_0; y) =& 4\delta_{\varepsilon=(1,1)} H_{w_l}^0(F_0; y) = 4\delta_{\varepsilon=(1,1)} \frac{1}{\abs{y_1 y_2}} \int_{\Re(\mu)=0} F_0(\mu) K^0_{w_l}(y, \mu) \specmu^0(\mu) d\mu,
\end{align*}
to which we may apply \eqref{eq:KwlWhitt}.
Thus for $y_1,y_2 > 0$, we have
\begin{align*}
	H_{w_l}^*(F_0; (y_1^2/4,y_2^2/4)) =& \frac{8\pi^5}{y_1 y_2} \int_{\Re(\mu)=0} F_0(\mu) W^{0*}(y, 2\mu) \sinmu^0(\mu) d\mu,
\end{align*}
recalling \eqref{eq:spec0}.

The double-angle formula gives
\[ 32\pi \cosmu^0(\mu) \sinmu^0(\mu) = \sinmu^0(2\mu), \]
so applying \eqref{eq:F0def} gives
\begin{align*}
	H_{w_l}^*(F_0; (y_1^2/4,y_2^2/4)) =& \frac{4}{y_1 y_2} \int_{\Re(\mu)=0} \int_{Y^+} t_1 t_2 f(t) W^{0*}(t,-2\mu) dt W^{0*}(y, 2\mu) \sinmu^0(2\mu) d\mu.
\end{align*}
Then Kontorovich-Lebedev inversion, Theorem \ref{thm:KontLeb}, gives
\[ H_{w_l}^*(F_0; (y_1^2/4,y_2^2/4)) = f(y). \]

\subsection{The other terms}
The trivial terms satisfy
\begin{align*}
	H_I^1(F_1; I) =& \int_{\Re(\mu)=0} F_0(\mu) \frac{1}{3} \sum_{w\in W_3} \frac{1}{2} \tan\frac{\pi}{2}(\mu^w_1-\mu^w_3)\tan\frac{\pi}{2}(\mu^w_2-\mu^w_3) \specmu^1(\mu^w) d\mu \\
	=& -H_I^0(F_0; I),
\end{align*}
which follows from the Weyl invariance of $F_0$ and the triple tangent identity
\[ \sum_{w \in W_3} \tan\frac{\pi}{2}(\mu^w_1-\mu^w_2) = -\prod_{i<j} \tan\frac{\pi}{2}(\mu_i-\mu_j). \]

Directly from the Weyl invariance of $F_0$ and equations \eqref{eq:Kw40Def}, \eqref{eq:Jw41Def}, \eqref{eq:Kw41Def}, \eqref{eq:spec0} and \eqref{eq:spec1} we see the $w_4$ weight functions satisfy
\begin{align*}
	H_{w_4}^0(F_0; y) =& \frac{3}{\abs{y_1}} \int_{\Re(\mu)=0} F_0(\mu) J_{w_4}(y, \mu) \sin\frac{\pi}{2}(\mu_1-\mu_2) \frac{\specmu^0(\mu)}{8\pi\prod_{i<j} \sin \frac{\pi}{2}\paren{\mu_i-\mu_j}} d\mu \\
	=& -H_{w_4}^1(F_1; y),
\end{align*}
and \eqref{eq:Kw5dDef} implies also
\[ H_{w_5}^0(F_0; y)+H_{w_5}^1(F_1; y) = 0. \]

This completes the proof of Theorem \ref{thm:ArithKuz}.

\section{The proof of Proposition \ref{prop:ArithKuzConv}}
\label{sect:ArithKuzConv}
Suppose $f,\sigma_1,\sigma_2,\theta$ and $N$ are as in the statement of the proposition and $\set{f_\ell}$ is any sequence of smooth, compactly supported functions on $(\R^+)^2$ tending to $f$ pointwise such that the derivatives $\partial_{y_1}^{j_1} \partial_{y_2}^{j_2} f_\ell$ also tend to $\partial_{y_1}^{j_1} \partial_{y_2}^{j_2} f$ pointwise.
We may assume $f$ and all $f_\ell$ satisfy conditions 1, 3 and 4 of the proposition with the same implied constants.
(The conditions may be expressed in terms of a collection of seminorms on an $L^1$-space to which we apply the usual density theory.)
The proof then consists of three applications of dominated convergence.

On the Kloosterman sum side, condition 1 of the proposition is sufficient since on $y=\paren{\frac{2\sqrt{m_1 n_2 c_2}}{c_1},\frac{2\sqrt{m_2 n_1 c_1}}{c_2}}$, we have
\[ (y_1 y_2)^{1+\epsilon} = \paren{\frac{16 m_1 m_2 n_1 n_2}{c_1 c_2}}^{\frac{1+\epsilon}{2}}, \]
which is sufficient for convergence of the sum of Kloosterman sums by \eqref{eq:StevensBd}.
(Note that $y_1 y_2 \to 0$ as $c_1,c_2 \to \infty$ so only the behavior as $y_1 y_2\to 0$ is relevant.)

On the spectral side, we proceed in two steps:
Let $\tilde{f}(y) = (y_1 y_2)^{-1} (\wtilde{\Delta_1})^N (y_1 y_2 f)$ and similarly for $\tilde{f}_\ell$, then with $\lambda_1(\mu)$ as in \eqref{eq:W0eigenfunc},
\begin{align}
\label{eq:ArithKuzConvStep1}
	\lim_{\ell\to\infty} F_0(f_\ell;\mu) = \lambda_1(-2\mu)^{-N} \lim_{\ell\to\infty} F_0(\tilde{f}_\ell;\mu) = \lambda_1(-2\mu)^{-N} F_0(\tilde{f};\mu) = F_0(f_\ell;\mu),
\end{align}
as we now explain.

The first and last equalities of \eqref{eq:ArithKuzConvStep1} hold because $\wtilde{\Delta_1}$ is self-adjoint with respect to the inner product
\[ \left<h_1,h_2\right>=\int_{Y^+} h_1(y) \wbar{h_2(y)} dy, \]
which can be seen directly through integration by parts or by noting that $\wtilde{\Delta_1}$ is a certain restriction of the Laplacian as in \cite[Theorem 6.1.6]{Gold01}.
(The elements of the operator algebra considered there are not typically self-adjoint, but the Laplacian certainly is.)
Condition 3 of the proposition assures that the boundary terms in the integration by parts are all zero since the Whittaker function has super-polynomial (in fact, exponential) decay as $y_i \to \infty$.
(Note that if some $y_i \to 0$ while $y_1 y_2 \gg 1$, then necessarily $y_{3-i} \to \infty$, so condition 3 is sufficient along the boundaries $y_i=0$ as well as $y_i=\infty$.)

The central equality of \eqref{eq:ArithKuzConvStep1} holds by dominated convergence using the bound of condition 4.
Here we need some bound on the Whittaker function, and we use \cite[Proposition 1]{Val01}:
\begin{thm}[Blomer]
\label{thm:WhittBd}
Suppose $\mu$ satisfies the unitaricity condition and the bound towards Ramanujan-Selberg \eqref{eq:muProps}.
Then for any $A > a_1 > \abs{a_2}+\theta$, $a_2 \in \R$, we have the bound
\[ W^{0*}(y,\mu) \ll_A (y_1 y_2)^{1-a_1} (y_1/y_2)^{a_2} (1+\norm{\mu})^{2a_1-\frac{1}{2}+\epsilon}. \]
\end{thm}

We split the integral defining $F_0(\tilde{f};\mu)$ or $F_0(\tilde{f}_\ell;\mu)$ into four pieces along the curves $y_1 y_2 = 1$ and $y_1=y_2$.
Then we apply Theorem \ref{thm:WhittBd} on each piece with $a_1=\sigma_1+\sigma_2\pm \epsilon$ according to $y_1 y_2 > 1$ or $y_1 y_2 < 1$ and $a_2 = \pm(\sigma_2 +\epsilon)$ according to $y_1 < y_2$ or $y_1 > y_2$.
By condition 4 of the proposition, the resulting integrals converge absolutely, and we have the bound
\begin{align}
\label{eq:F0Bd}
	\lambda_1(-2\mu)^{-N} F_0(\tilde{f};\mu) \ll \lambda_1(-2\mu)^{-\sigma_1-\sigma_2-\frac{9}{4}-10\epsilon} (1+\norm{\mu})^{2\sigma_1+2\sigma_2-\frac{1}{2}+\epsilon} \ll (1+\norm{\mu})^{-5-\epsilon},
\end{align}
and the same for $\lambda_1(-2\mu)^{-N} F_0(\tilde{f}_\ell;\mu)$, with the same implied constants.
From the Weyl law \cite{WeylI}*{Theorem 1} and its spherical counterpart \cite[Theorem 1]{Val01} (which also dominate the Weyl laws for the continuous spectra), a bound of $(1+\norm{\mu})^{-5-\epsilon}$ is sufficient for dominated convergence on the spectral expansion, so we are done.

\section{The poles, residues and decay of $G^0$}
\label{sect:Residues}
Define
\begin{align*}
	\wtilde{F}_0(u,\mu) =& \frac{4}{\pi^4} \cosmu^0(\mu) G^0(2u,-2\mu), 
\end{align*}
and by analogy,
\begin{align*}
	\wtilde{F}_1(u,\mu) =& \frac{1}{2} \tan\frac{\pi}{2}(\mu_1-\mu_3)\tan\frac{\pi}{2}(\mu_2-\mu_3) \wtilde{F}_0(u,\mu).
\end{align*}
These functions are connected to the integral transforms $F_d(f;\mu)$ through Plancherel's theorem for the Mellin transform.

Since the tangents are particularly non-threatening, we analyze $\wtilde{F}_0$ and leave the other to the reader.
First, we note that Stirling's formula implies $\wtilde{F}_0(u,\mu)$, for $\Re(\mu)$ and $\Re(u)$ in fixed compact sets, has no exponential growth and in fact decays exponentially in $\Im(\mu)$ or $\Im(u)$ unless
\begin{align}
\label{eq:G0noExpDecay}
	\Im(\mu_1) \le \Im(u_2) \le \Im(\mu_2) \le \Im(-u_1) \le \Im(\mu_3),
\end{align}
up to permutation of the coordinates of $\mu$ or the coordinates of $(-u_1,u_2)$.
This is somewhat unpleasant to check, so we refer the reader to \cite{Val01}*{eqs. (2.19)-(2.21)}.
In particular, $\wtilde{F}_0(u,(it,it,-2it))$ decays exponentially as $t \to \pm\infty$.
Furthermore, $\mu_1+\mu_2+\mu_3=0$ implies
\[ \prod_{i=1}^3 (1+\abs{\mu_i}) \gg (1+\norm{\mu})^2, \]
so for, say, $\norm{\Im(u)} \ll (1+\norm{\mu})^\epsilon$, $\abs{\Re(\mu)} < \frac{1}{2}$ and $\Re(u)$ in some fixed compact set with $\Re(u_1+u_2) < 0$, we have
\[ \wtilde{F}_0(u,\mu) \ll (1+\norm{\mu})^{2\Re(u_1+u_2)+\epsilon}, \]
even when one of the coordinates of $\mu$ is small (say $\Im(\mu_i) \ll (1+\norm{\mu})^\epsilon$).

Next we compute the residues of $\wtilde{F}_0(\tilde{s},\mu)$:
When the coordinates of $\mu$ are distinct modulo $\Z$, we have
\begin{gather}
\label{eq:tildeF0ResGenMu}
	\res_{u_1=-\mu_1-\ell} \wtilde{F}_0(u,\mu) = \frac{4}{\pi^4} \cosmu^0(\mu) \frac{(-1)^\ell}{\ell!} \frac{\Gamma(\mu_2-\mu_1-\ell)\Gamma(\mu_3-\mu_1-\ell)}{\Gamma(u_2-\mu_1-\ell)} \prod_{i=1}^3 \Gamma(u_2-\mu_i), \\
\begin{aligned}
\label{eq:tildeF0ResGenMuDouble}
	\res_{u_1=-\mu_1-\ell_1} \res_{u_2=\mu_3-\ell_2} \wtilde{F}_0(u,\mu) =& \frac{4}{\pi^4} \cosmu^0(\mu) \frac{(-1)^{\ell_1+\ell_2}}{\ell_1!\,\ell_2!}\frac{\Gamma\paren{\mu_3-\mu_1-\ell_1} \Gamma\paren{\mu_3-\mu_1-\ell_2}}{\Gamma\paren{\mu_3-\mu_1-\ell_1-\ell_2}} \\
	& \times  \Gamma\paren{\mu_2-\mu_1-\ell_1} \Gamma\paren{\mu_3-\mu_2-\ell_2}.
\end{aligned}
\end{gather}
Notice that the exponential parts of the second-order residue \eqref{eq:tildeF0ResGenMuDouble} cancel.

When $\mu=(it,it,-2it)$ with $0\ne t\in\R$, the residue at $u_1=2it-\ell$ follows the form \eqref{eq:tildeF0ResGenMu}, and the residue at the double pole becomes
\begin{equation}
\label{eq:tildeF0ResDegenMu}
\begin{aligned}
	\res_{u_1=-it-\ell} \wtilde{F}_0(u,\mu) =& \frac{4}{\pi^4} \cosmu^0(\mu) \paren{\prod_{i=1}^3 \Gamma(u_2-\mu_i)} \frac{1}{(\ell!)!} \frac{\Gamma(-3it-\ell)}{\Gamma(u_2-it-\ell)} \\
	& \qquad \times \paren{2H_\ell-2\gamma+\psi(-3it-\ell)-\psi(u_2-it-\ell)},
\end{aligned}
\end{equation}
where $H_\ell = \sum_{j=1}^\ell \frac{1}{j}$ is the $\ell$-th harmonic number, $\gamma$ is the Euler-Mascheroni constant and $\psi(z)=\frac{\Gamma'(z)}{\Gamma(z)}$ is the digamma function.
Since computing residues at double poles requires differentiation, we do not include the second-order residues here.

Lastly, $\wtilde{F}_0(u,0)$ has a triple pole at $u_1=-\ell$, $\ell\in\N_0$; the residue is easily computed with a computer algebra package and is somewhat complicated, so we do not include it here.
As pointed out in Section \ref{sect:KloosZeta}, the case $\mu=0$ is not relevant to the spectral expansion.

\section{The proof of Theorem \ref{thm:SmoothSums}}
\label{sect:SmoothSums}
We apply Theorem \ref{thm:ArithKuz} to the function
\[ f_X(y) = f\paren{\frac{X_1 y_1^2}{4 m_1 n_2}, \frac{X_2 y_2^2}{4 m_2 n_1}}. \]
Using the Mellin-Barnes integral \eqref{eq:SphWhitt} (i.e. by Plancherel's theorem for the Mellin transform), we write $F_0(f_X; \mu)$ in the form
\begin{align}
\label{eq:SmoothSumsF0MB}
	F_0(f_X; \mu) =& \int_{\Re(s)=(\theta+\epsilon,\theta+\epsilon)} \paren{\frac{X_1}{4\pi^2 m_1 n_2}}^{s_1} \paren{\frac{X_2}{4\pi^2 m_2 n_1}}^{s_2} \hat{f}\paren{-s} \wtilde{F}_0(s,\mu) \frac{ds_1 \, ds_2}{(2\pi i)^2},
\end{align}
where $\hat{f}(s)$ is the usual Mellin transform of $f$.
Now $\hat{f}(s)$ has super-polynomial decay in $\Im(s)$, and, as described in Section \ref{sect:Residues}, $\wtilde{F}_0(s,\mu)$ has exponential decay in $\Im(\mu)$ unless some $\abs{\Im(\mu_i)} \ll \norm{\Im(s)}$.
Therefore,
\[ F_0(\mu) \ll_{m,n,f} (X_1 X_2 (1+\norm{\mu}))^{-100} \]
unless some $\abs{\Im(\mu_i)} \ll (X_1 X_2)^\epsilon$.

From the Weyl law \cite{WeylI}*{Theorem 1} and its spherical counterpart \cite[Theorem 1]{Val01}, the tempered cusp forms and continuous spectrum of Theorem \ref{thm:ArithKuz} with $\norm{\mu} \ll (X_1 X_2)^\epsilon$ contribute at most $(X_1 X_2)^\epsilon$ to the bound in Theorem \ref{thm:SmoothSums}, which is observable from the trivial bound on \eqref{eq:SmoothSumsF0MB} with the contours at $\Re(s) = (\epsilon,\epsilon)$.
If there exist any non-tempered (i.e. complementary series) cusp forms with spectral parameters, say, $\mu$ a permutation of $(-x+it,-2it,x+it)$, $0<x<\theta$ and $t \ll (X_1 X_2)^\epsilon$ their contribution is at most $(X_1 X_2)^{\theta+\epsilon}$.

For the troublesome forms which are near self-dual forms, say with one spectral parameter $\abs{\Im(\mu_i)} \ll (X_1 X_2)^\epsilon$ and $\norm{\mu} \gg (X_1 X_2)^{2\epsilon}$ (which rules out the non-tempered forms), we suppose WLOG $i=1$ and $X_1 > X_2$ (so, in particular $X_1 > 1$) and shift $\Re(s_1)$ highly negative (saving arbitrarily many powers of $X_1 (1+\norm{\mu})^2$), giving
\begin{equation}
\label{eq:SmoothSumsF0Asymp}
\begin{aligned}
	F_0(f_X; \mu) \approx& \sum_{w\in W_3} \paren{\frac{X_1}{4\pi^2 m_1 n_2}}^{-\mu_1^w} \int_{\Re(s_2)=\epsilon} \paren{\frac{X_2}{4\pi^2 m_2 n_1}}^{s_2} \hat{f}\paren{\mu_1^w,s_2} \res_{s_1=-\mu_1^w} \wtilde{F}_0(s,\mu) \frac{ds_1 \, ds_2}{(2\pi i)^2} \\
	&+\text{lower-order terms}.
\end{aligned}
\end{equation}
All of the terms have super-polynomial decay in $\mu_2\asymp-\mu_3$ (coming from $\hat{f}$) except $w=I$, and on that term, we shift $\Re(s_2) \mapsto -1-\epsilon$.
The residues again have super-polynomial decay in $\mu_2$, and the integral on the shifted contour is $\ll X_2^{-1-\epsilon} \norm{\mu}^{-4-\epsilon}$.
This is sufficient for this part of the spectral expansion to converge, as the Weyl laws \cite{Val01}*{Theorem 1} and \cite{WeylI}*{Theorem 1} imply there are at most $T^{4+\epsilon}(X_1 X_2)^\epsilon$ forms with $\mu_1 \ll (X_1 X_2)^\epsilon$ and $\mu_2 \asymp -\mu_3 \asymp T$.

\section{The proof of Proposition \ref{prop:KloosZetaTechnical}}
\label{sect:KloosZetaTechnical}

Proposition \ref{prop:ArithKuzConv} applies for $f=f_s$ with $\sigma_1=\theta+\epsilon$, $\sigma_2=\epsilon$, $N=3$ on the region $2\Re(\tilde{s}_1),2\Re(\tilde{s}_2) > 1$, which proves part 1.
The remaining parts concern
\begin{align}
\label{eq:F0smuMB}
	F_0(s,\mu) =& \int_{\Re(u_1,u_2)=(\epsilon,\epsilon)} \Gamma\paren{u_1} \Gamma\paren{u_2} \wtilde{F}_0(\tilde{s}-u,\mu) \frac{du_1 \, du_2}{(2\pi i)^2},
\end{align}
after applying \eqref{eq:SphWhitt} as in the previous section.

When the components of $\mu$ are distinct modulo $\Z$, we have
\begin{equation}
\label{eq:F0muiDistinct}
\begin{aligned}
	F_0(s,\mu) =& \sum_{w\in W} \sum_{j_1,j_2 \ge 0} \Gamma\paren{j_1+\tilde{s}_1+\mu^w_1} \Gamma\paren{j_2+\tilde{s}_2-\mu^w_3} \res_{u_1=-\mu_1^w-j_1} \res_{u_2=\mu_3^w-j_2} \wtilde{F}_0(u,\mu),
\end{aligned}
\end{equation}
by shifting the $u$ contours to $\infty$, as we may.
From \eqref{eq:tildeF0ResGenMuDouble}, this series representation converges rapidly away from the poles and clearly defines a meromorphic function of all $\tilde{s}$ and $\mu$, proving part 2, and part 4 is clear.
Note that the product of $\cosmu^0(\mu)$ and the gamma factors of the summand has exponential decay for $\Im(\mu)$ large compared to $\Im(\tilde{s})$ and away from the poles, which is the content of part 3, so this is proved for $\mu$ whose components are distinct modulo $\Z$.
There is a minor detail that when, say, $\mu_1\asymp\mu_2$ is large compared to $\mu_1-\mu_2$, the polynomial part of the summand is increasing in $j_1$ like $\abs{\mu_1}^{j_1}$, but to defeat the exponential decay, we would need $j \sim \frac{\abs{\mu_1}}{\log\abs{\mu_1}}$, at which point the $j_1!$ in the denominator (as well as the factor $\Gamma\paren{\mu_2-\mu_1-j_1}$) overwhelms any possible growth.

When $\mu=(it,it,-2it)$ with $t \ne 0$, we have
\begingroup\allowdisplaybreaks\begin{gather}
\label{eq:F0mutt2t}
	F_0(s,\mu) = \frac{2}{\pi^4} \cosmu^0(\mu) \sum_{i=1}^3 F_{0,i}(s,t), \\
\begin{aligned}
	F_{0,1}(s,t) =& \sum_{j_1,j_2\ge 0} \frac{(-1)^{j_1+j_2}(j_1+j_2)!}{(j_1!\,j_2!)^2} \Gamma(-3it-j_1)\Gamma(3it-j_2)\Gamma(\tilde{s}_1+it+j_1)\Gamma(\tilde{s}_2-it+j_2) \nonumber \\
	& \qquad \times\bigl(2H_{j_2}-2H_{j_1+j_2}+2\psi(j_1+1)+\psi(-3it-j_1)+\psi(3it-j_2) \nonumber \\
	& \qquad -\psi(\tilde{s}_1+it+j_1)-\psi(\tilde{s}_2-it+j_2)\bigr) \nonumber
\end{aligned}\\
\begin{aligned}
	F_{0,2}(s,t) =& \sum_{j_1,j_2\ge 0} \frac{(-1)^{j_2}}{(j_1!)^2\,j_2!} \frac{\Gamma(-3it-j_1)\Gamma(-3it-j_2)^2}{\Gamma(-3it-j_1-j_2)} \Gamma(\tilde{s}_1+it+j_1)\Gamma(\tilde{s}_2+2it+j_2) \nonumber \\
	& \qquad \times\bigl(2\psi(j_1+1)+\psi(-3it-j_1) -\psi(\tilde{s}_1+it+j_1)-\psi(-3it-j_1-j_2)\bigr) \nonumber \\
	F_{0,3}(s,t) =& F_{0,2}((s_2,s_1),-t), \nonumber
\end{aligned}
\end{gather}\endgroup
and it is easy to see that $F_0(s,\mu)$ decays exponentially as $t \to \pm\infty$, which completes the proof of part 3.

The equations \eqref{eq:F0muiDistinct} and \eqref{eq:F0mutt2t} were necessary to show the rapid decay in $\mu$ for part 3 (for all $s\in\C^2$ away from the poles), but to see the poles with their locations, orders and residues, it is perhaps easier to write
\begin{equation}
\label{eq:F0poles}
\begin{aligned}
	F_0(s,\mu) =& \sum_{j_1=0}^{T_1} \sum_{j_2=0}^{T_2} \frac{(-1)^{j_1+j_2}}{j_1!\,j_2!} \wtilde{F}_0(\tilde{s}+j,\mu) \\
	&+\sum_{j_2=0}^{T_2} \frac{(-1)^{j_2}}{j_2!} \int_{\Re(u_1)=-T_1-\frac{1}{2}} \Gamma(u_1) \wtilde{F}_0(\tilde{s}+(-u_1,j_2),\mu) \frac{du_1}{2\pi i} \nonumber \\
	&+\sum_{j_1=0}^{T_1} \frac{(-1)^{j_1}}{j_1!} \int_{\Re(u_2)=-T_2-\frac{1}{2}} \Gamma(u_2) \wtilde{F}_0(\tilde{s}+(j_1,-u_2),\mu) \frac{du_2}{2\pi i} \nonumber \\
	&+ \int_{\Re(u)=(-T_1-\frac{1}{2},-T_2-\frac{1}{2})} \Gamma(u_1) \Gamma(u_2) \wtilde{F}_0(\tilde{s}-u,\mu) \frac{du_1\,du_2}{(2\pi i)^2}, \nonumber
\end{aligned}
\end{equation}
assuming $-\Re(\tilde{s}_i) < T_i \in \Z$.
Now the poles are easy to compute from the poles of $\wtilde{F}_0(\tilde{s},\mu)$; for instance, if $\ell\in\N_0$, then
\begin{align}
\label{eq:F0fromTildeF0}
	\res_{\tilde{s}_1=-\mu_j-\ell} F_0(s,\mu) =& \sum_{j_1=0}^{\ell} \sum_{j_2=0}^{T_2} \frac{(-1)^{j_1+j_2}}{j_1!\,j_2!} \res_{\tilde{s}_1=-\mu_j-\ell} \wtilde{F}_0(\tilde{s}+j,\mu) \\
	&+\sum_{j_1=0}^{\ell} \frac{(-1)^{j_1}}{j_1!} \int_{\Re(u_2)=-T_2-\frac{1}{2}} \Gamma(u_2) \res_{\tilde{s}_1=-\mu_j-\ell} \wtilde{F}_0(\tilde{s}+(j_1,-u_2),\mu) \frac{du_2}{2\pi i}, \nonumber
\end{align}
to which we apply \eqref{eq:tildeF0ResGenMu} or \eqref{eq:tildeF0ResDegenMu}.

\bibliographystyle{amsplain}

\bibliography{HigherWeight}

\end{document}